\documentclass[12pt,a4paper,english]{article}
\usepackage[T1]{fontenc}
\usepackage[latin9]{inputenc}
\usepackage{amsmath}
\usepackage{amssymb}
\usepackage{esint}

\makeatletter


\makeatother

\usepackage{babel}

\makeatother

\usepackage{babel}

\begin{document}

\title{Multiplicative white noise functionals and the Krylov-Veretennikov
expansion for coalescing stochastic flows}

\author{Andrey A. Dorogovtsev%
\thanks{This article was partially supported by the State fund for fundamental
researches of Ukraine and The Russian foundation for basic research,
grant F40.1/023.%
}}
\maketitle
\begin{abstract}
In this article we consider multiplicative operator-valued white noise
functionals related to a stochastic flow. A generalization of the
Krylov-Veretennikov expansion is presented. An analog of such expansion
for the Arratia flow is derived.\end{abstract}
\begin{description}
\item [{{\large Introduction.}}]~{\large \par}
\end{description}
{\large In this article we present the form of the kernels in the
Itô-Wiener expansion for functionals from a dynamical system driven
by an additive Gaussian white noise. The most known example of such
expansion is the Krylov-Veretennikov representation {[}1{]}: \[
f(y(t))=T_{t}f(u)+\sum_{k=1}^{\infty}\int_{\Delta_{k}(0;t)}T_{\tau_{1}}b\partial T_{\tau_{2}-\tau_{1}}\ldots\]
 \[
\ldots b\partial T_{t-\tau_{k}}f(u)dw(\tau_{1})\ldots dw(\tau_{k}).\]
 where $f$ is a bounded measurable function, $y$ is a solution of
SDE \[
dy(t)=a(y(t))dt+b(y(t))dw(t)\]
 with smooth and nondegenerate coefficients, and $\{T_{t};\; t\geq0\}$
is the semigroup of operators related to SDE and $\partial$ is the
symbol of differentiation. }{\large \par}

{\large A family of substitution operators of SDE's solution into
a function can be treated as a multiplicative Gaussian white noise
functional. In the first section of this article we consider a family
$\{G_{s,t};\:0\leq s\leq t<+\infty\}$ of strong random operators
in the Hilbert space which is an operator-valued multiplicative functional
from the Gaussian white noise. It turns out that the precise form
of the kernels in the Itô-Wiener expansion can be found for a wide
class of operator-valued multiplicative functionals using some simple
algebraic relations. The obtained formula covered the Krylov-Veretennikov
case and gives a representation for different objects such as Brownian
motion in Lie group etc.. }{\large \par}

{\large The representation obtained in the first section may be useful
in studing the properties of a dynamical system with an additive Gaussian
white noise. On the other hand, there exist cases when a dynamical
system is obtained as a limit in a certain sense of systems driven
by the Gaussian white noise. A limiting system could be highly irregular
{[}4, 12, 13{]}. One example of such system is the Arratia flow {[}12{]}
of coalescing Brownian particles on the real line. The trajectories
of individual particles in this flow are Brownian motions, but the
whole flow cannot be built from the Gaussian noise in a regular way
{[}14{]}. Nevertheless, it is possible to construct the $n-$point
motion of the Arratia flow from the pieces of the trajectories of
$n$ independent Wiener processes. Correspondingly a function from
the $n-$point motion of the Arratia flow has an Itô-Wiener expansion
based on the initial Wiener processes. This expansion depends on the
way of construction (coalescing description). We present such expansion
in terms of an infinite family of expectation operators related to
all manner of coalescence of the trajectories in the Arratia flow.
To do this we first obtain an analog of the Krylov-Veretennikov expansion
for the Wiener process stopped at zero.}{\large \par}

{\large The article is divided onto three parts. The first section
is devoted to multiplicative operator-valued functionals from Gaussian
white noise. The second part contains the definition and necessary
facts about the Arratia flow. In the last section we present a family
of Krylov-Veretennikov expansions for the $n-$point motion of the
Arratia flow.}{\large \par}

\section{Multiplicative white noise functionals}

{\large In this part we present the Itô-Wiener expansion for the semigroup
of strong random linear operators in Hilbert space. Such operators
in the space of functions can be generated by the flow of solutions
to a stochastic differential equation. In this case our expansion
turns into the well-known Krylov--Veretennikov representation {[}1{]}.
In the case when these operators have a different origine, we obtain
a new representation for the semigroup.}{\large \par}

{\large Let us start with the definition and examples of strong random
operators in the Hilbert space. Let $H$ denote a separable real Hilbert
space with norm $\|\cdot\|$ and inner product $(\cdot,\cdot).$ As
usual $(\Omega,\mathcal{F},P)$ denotes a complete probability space.}{\large \par}

\begin{description}
\item [{{\large Definition}}] \textbf{\large 1.1.}{\large{} A strong linear
random operator in $H$ is a continuous linear map from $H$ to $L_{2}(\Omega,P,H).$ }{\large \par}
\end{description}
\noindent \textbf{\textit{\large Remark.}}{\large{} The notion of strong
random operator was introduced by A.V.Skorokhod {[}2{]}. In his definition
Skorokhod used the convergence in probability, rather than convergence
in the square mean. }{\large \par}

{\large Consider some typical examples of strong random operators. }{\large \par}

\textbf{\large Example 1.1.}{\large{} Let $H$ be $l_{2}$ with the
usual inner product and $\{\xi_{n};\ n\geq1\}$ be an i.i.d. sequence
with finite second moment. Then the map \[
l_{2}\ni x=(x_{n})_{n\geq1}\mapsto Ax=(\xi_{n}x_{n})_{n\geq1}\]
 is a strong random operator. Really, \[
E\|Ax\|^{2}=\sum_{n=1}^{\infty}x_{n}^{2}E\xi_{1}^{2}\]
 and the linearity is obvious. Note, that pathwise the operator $A$
can be not well-defined. For example, if $\{\xi_{n}:n\geq1\}$ have
the standard normal distribution, then with probability one \[
\sup_{n\geq1}|\xi_{n}|=+\infty.\]
 An interesting set of examples of strong random operators can be
found in the theory of stochastic flows. Let us recall the definition
of a stochastic flow on $\mathbb{R}$ {[}3{]}.}{\large \par}

\textbf{\large Definition 1.2. }{\large A family $\{\phi_{s,t};\:0\leq s\leq t\}$
of random maps of $\mathbb{R}$ to itself is referred to as a stochastic
flow if the following conditions hold}{\large \par}

\begin{enumerate}
\item {\large For any $0\leq s_{1}\leq s_{2}\leq\ldots s_{n}<\infty:\;\phi_{s_{1},s_{2}}$
,$\ldots$ ,$\phi_{s_{n-1},s_{n}}$ are independent.}{\large \par}
\item {\large For any $s,t,r\geq0:\;$$\phi_{s,t}$ and $\phi_{s+r,t+r}$
are equidistributed.}{\large \par}
\item {\large For any $r\leq s\leq t$ and $u\in\mathbb{R}:\;$ $\phi_{r,s}\phi_{s,t}(u)=\phi_{r,t}(u)$,
$\phi_{r,r}$ is an identity map.}{\large \par}
\item {\large For any $u\in\mathbb{R}:\;$$\phi_{0,t}(u)\mapsto u$ in probability
when $t\mapsto0.$ }{\large \par}
\end{enumerate}
{\large Stochastic flows arise as solutions to stochastic differential
equations with smooth coefficients. Namely, if $\phi_{s,t}(u)$ is
a solution to the stochastic differential equation\[
dy(t)=a(y(t))dt+b(y(t))dw(t)\qquad(1.1)\]
 starting at the point $u$ in time $s$ and considered in time $t$,
then under smoothness conditions on the coefficients $a$ and $b$
the family $\{\phi_{s,t}\}$ will satisfy the conditions of Definition
1.2 {[}3{]}. Another example of a stochastic flow is the Harris flow
consisting of Brownian particles {[}4{]}. In this flow $\phi_{0,t}(u)$
for every $u\in\mathbb{R}$ is a Brownian martingale with respect
to a common filtration and \[
d<\phi_{0,t}(u_{1}),\:\phi_{0,t}(u_{2})>=\Gamma(\phi_{0,t}(u_{1})-\phi_{0,t}(u_{2}))dt\]
 for some positive definite function $\Gamma$ with $\Gamma(0)=1.$}{\large \par}

 {\large For a given stochastic flow one can try to construct corresponding
family of strong random operators as follows.}{\large \par}

\textbf{\large Example 1.2.}{\large{} Let $H=L_{2}(\mathbb{R})$. Define
\[
G_{s,t}f(u)=f(\phi_{s,t}(u)).\]
 Let us check that in the both cases mentioned above $G_{s,t}$ satisfies
Definition 1.1. For the Harris flow we have \[
E\int_{\mathbb{R}}f(\phi_{s,t}(u))^{2}du=\int_{\mathbb{R}}\int_{\mathbb{R}}f(v)^{2}p_{1}(u-v)dudv=\int_{\mathbb{R}}f(v)^{2}dv.\]
 Here $p_{r}$ denotes the Gaussian density with zero mean and variance
$r.$}{\large \par}

{\large To get an estimation for the flow generated by a stochastic
differential equation let us suppose that the coefficients $a$ and
$b$ are bounded Lipschitz functions and $b$ is separated from zero.
Under such conditions $\phi_{s,t}(u)$ has a density, which can be
estimated from above by a Gaussian density {[}5{]}. Consequently we
will have the inequality $E\int_{\mathbb{R}}f(\phi_{s,t}(u))^{2}du\leq\int_{\mathbb{R}}f(v)^{2}dv$.}{\large \par}

{\large As it was shown in Example 1.1, a strong random operator in
general is not a family of bounded linear operators in $H$ indexed
by the points of probability space. Despite of this the superposition
of such operators can be properly defined (see {[}6{]} for detailed
construction in case of dependent nonlinear operators via Wick product).
Here we will consider only the case when strong random operators $A$
and $B$ are independent. In this case both $A$ and $B$ have measurable
modifications and one can define for $u\in H,\:\omega\in\Omega$ \[
AB(u,\omega):=A(B(u,\omega),\omega)\]
 and prove that the value $AB(u)$ does not depend on the choice of
modifications. Note, that the operators from the previous example
satisfy the semigroup property, and that for the flow generated by
a stochastic differential equation these operators are measurable
with respect to increments of the Wiener process. In this section
we will consider a general situation of this kind and study the structure
of the semigroup of strong random operators measurable with respect
to a Gaussian white noise.}{\large \par}

{\large Let's start with a description of the noise. Let $H_{0}$
be a separable real Gilbert space. Define a new Hilbert space $\tilde{H}=H_{0}\otimes L_{2}([0;\:+\infty]),$
where an inner product is defined by the formula\[
\tilde{H}\ni f,g\mapsto<f,\: g>=\int_{0}^{\infty}(f(t),\: g(t))_{0}dt.\]
 }{\large \par}

\textbf{\large Definition 1.3.}{\large{} A family of jointly Gaussian
random variables $\{<\xi,h>;\: h\in\tilde{H}\}$ is referred to as
a Gaussian white noise in $\tilde{H}$ if it is linear with respect
to $h\in\tilde{H}$ and for every $h$, $<\xi,h>$ has zero mean and
variance $\|h\|^{2}$.}{\large \par}

{\large Let $\tilde{H}_{s,t}$ be the product $H_{0}\otimes L_{2}([s;\: t]),$
which can be naturally considered as a subspace of $\tilde{H}.$ Define
the $\sigma$-fields $\mathcal{F}_{s,t}=\sigma\{<\xi,h>;\: h\in\tilde{H}_{s,t}\},\:0\leq s\leq t<+\infty.$ }{\large \par}

\textbf{\large Definition 1.4.}{\large{} A family $\{G_{s,t};\:0\leq s\leq t<+\infty\}$
of strong random operators in $H$ is refereed to as a multiplicative
functional from $\xi$ if the following conditions hold:}{\large \par}

{\large 1) $G_{s,t}$ is measurable with respect to $\mathcal{F}_{s,t},$}{\large \par}

{\large 2) $G_{s,s}$ is an identity operator for every $s,$}{\large \par}

{\large 3) $G_{s_{1},s_{3}}=G_{s_{2},s_{3}}G_{s_{1},s_{2}}$ for $s_{1}\leq s_{2}\leq s_{3}.$ }{\large \par}

\textbf{\large Remark. }{\large Taking an orthonormal basis $\{e_{n}\}$
in $H_{0}$ one can replace $\xi$ by a sequence of independent Wiener
processes $\{w_{n}(t)=<e_{n}\otimes1_{[0;\: t]};\:\xi>;\; t\geq0\}$.
We use the notation $\xi$ in order to simplify notations and consider
simultaneously both cases of finite and infinite number of the processes
$\{w_{n}\}.$}{\large \par}

\textbf{\large Example 1.3.}{\large{} Let us define $x(u,s,t)$ as
a solution to Cauchy problem for (1.1) which starts from the point
$u$ at the moment $s.$ Using the flow property one can easily verify
that the family of operators $\{G_{s,t}f(u)=f(x(u,s,t))\}$ in $L_{2}(\mathbb{R})$
is a multiplicative functional from the Gaussian white noise $\dot{w}$
in $L_{2}([0;+\infty]).$}{\large \par}

{\large Now we are going to introduce the notion of a homogeneous
multiplicative functional. Let us recall, that every square integrable
random variable $\alpha$ measurable with respect to $\xi$ can be
uniquely expressed as a series of multiple Wiener integrals {[}7{]}
\[
\alpha=E\alpha+\sum_{k=1}^{\infty}\int_{\Delta_{k}(0;+\infty)}a_{k}(\tau_{1},\ldots,\tau_{k})\xi(d\tau_{1})\ldots\xi(d\tau_{k}),\]
 where \[
\Delta_{k}(s;t)=\{(\tau_{1},\ldots,\tau_{k}):\ s\leq\tau_{1}\leq\ldots\leq\tau_{k}\leq t\},\]
 \[
a_{k}\in L_{2}(\Delta_{k}(0;+\infty),\: H_{0}^{\otimes k}),\: k\geq1.\]
 Here in the multiple integrals we consider the white noise $\xi$
as Gaussian $H_{0}-$valued random measure on $[0;\:+\infty).$ In
the terms of the mentioned above orthonormal basis $\{e_{n}\}$ in
$H_{0}$ and the sequence of the independent Wiener processes $\{w_{n}\}$
one can rewrite the above multiple integrals as \[
\int_{\Delta_{k}(0;+\infty)}a_{k}(\tau_{1},\ldots,\tau_{k})\xi(d\tau_{1})\ldots\xi(d\tau_{k})=\]
\[
\sum_{n_{1},...,n_{k}}\int_{\Delta_{k}(0;+\infty)}a_{k}(\tau_{1},\ldots,\tau_{k})(e_{n_{1}},...,e_{n_{k}})dw_{n_{1}}(\tau_{1})...dw_{n_{k}}(\tau_{k}).\]
 Define the shift of $\alpha$ for $r\geq0$ as follows \[
\theta_{r}\alpha=E\alpha+\sum_{k=1}^{\infty}\int_{\Delta_{k}(r;+\infty)}a_{k}(\tau_{1}-r,\ldots,\tau_{k}-r)\xi(d\tau_{1})\ldots\xi(d\tau_{k}).\]
 }{\large \par}

\textbf{\large Definition 1.5.}{\large{} A multiplicative functional
$\{G_{s,t}\}$ is homogeneous if for every $s\leq t$ and $r\geq0$
\[
\theta_{r}G_{s,t}=G_{s+r,t+r}.\]
 }{\large \par}

{\large Note, that the family $\{G_{s,t}\}$ from Example 1.3 is a
homogeneous functional. From now on, we will consider only homogeneous
multiplicative functionals from $\xi.$ For a homogeneous functional
$\{G_{s,t}\}$ one can define the expectation operators \[
T_{t}u=EG_{0,t}u,\ u\in H,\: t\geq0.\]
 Since the family $\{G_{s,t}\}$ is homogeneous, then $\{T_{t}\}$
is the semigroup of bounded operators in $H$. Under the well-known
conditions the semigroup $\{T_{t}\}$ can be described by its generator.
However the family $\{G_{s,t}\}$ cannot be recovered from this semigroup.
The following simple example shows this. }{\large \par}

\textbf{\large Example 1.4.}{\large{} Define $\{G_{s,t}^{1}\}$ and
$\{G_{s,t}^{2}\}$ in the space $L_{2}(\mathbb{R})$ as follows \[
G_{s,t}^{1}f(u)=T_{t-s}f(u),\]
 where $\{T_{t}\}$ is the heat semigroup, and \[
G_{s,t}^{2}f(u)=f(u+w(t)-w(s)),\]
 where $w$ is a standard Wiener processes. It is evident, that \[
EG_{s,t}^{2}f(u)=T_{t-s}f(u)=EG_{s,t}^{1}f(u).\]
 To recover multiplicative functional uniquely we have to add some
information to $\{T_{t}\}.$ It can be done in the following way.
For $f\in H$ define an operator which acts from $H_{0}$ to $H$
by the rule \[
A(f)(h)\doteq\lim_{t\to0+}\frac{1}{t}EG_{0,t}f(\xi,\: h\otimes1_{[0;\: t]}).\qquad(1.2)\]
 }{\large \par}

\textbf{\large Example 1.5.}{\large{} Let the family $\{G_{s,t}\}$
be defined as in Example 1.3. Now $H=L_{2}(\mathbb{R})$ and the noise
$\xi$ is defined on $L_{2}([0;\:+\infty)$ as $\dot{w}.$ Then for
$f\in L_{2}(\mathbb{R})$ (now $H_{0}=\mathbb{R}$ and it makes sense
only to take $h=1$) \[
A(f)(u)=\lim_{t\to0+}\frac{1}{t}Ef(x(u,t))w(t).\]
 Suppose that $f$ has two bounded continuous derivatives. Then using
Itô's formula one can get \[
Ef(x(u,t))w(t)=\int_{0}^{t}Ef'(x(u,s))\varphi(x(u,s))ds,\]
 and \[
\frac{1}{t}Ef(x(\bullet,t))w(t)\overset{L_{2}(\mathbb{R})}{\rightarrow}f'(\bullet)b(\bullet),\: t\to0+.\]
 Consequently, for {}``good'' functions \[
Af=bf'.\]
 }{\large \par}

\textbf{\large Definition 1.6.}{\large{} An element $u$ of $H$ belongs
to the domain of definition $D(A)$ of $A$ if the limit (1.2) exists
for every $h\in H_{0}$ and defines a Hilbert--Schmidt operator $A(u):H_{0}\to H.$
The operator $A$ is refereed to as the random generator of $\{G_{s,t}\}.$ }{\large \par}

{\large Now we can formulate the main statement of this section, which
describes the structure of homogeneous multiplicative functionals
from $\xi.$ }{\large \par}

\textbf{\large Theorem 1.1.}{\large{} Suppose, that for every $t>0$,
$T_{t}(H)\subset D(A)$ and the kernels of the Itô-Wiener expansion
for $G_{0,t}$ are continuous with respect to time variables. Then
$G_{0,t}$ has the following representation \[
G_{0,t}(u)=T_{t}u+\]
\[
\sum_{k=1}^{\infty}\int_{\Delta_{k}(0;t)}T_{t-\tau_{k}}AT_{\tau_{k}-\tau_{k-1}}\ldots AT_{\tau_{1}}ud\xi(\tau_{1})\ldots d\xi(\tau_{k}).\qquad(1.3)\]
 }{\large \par}

\textbf{\large Proof.}{\large{} Let us denote the kernels of the Itô-Wiener
expansion for $G_{0,t}(u)$ as $\{a_{k}^{t}(u,\tau_{1},\ldots,\tau_{k});\: k\geq0\}.$
Since \[
a_{0}^{t}(u)=EG_{0,t}(u),\]
 then \[
a_{0}^{t}(u)=T_{t}u.\]
 Since \[
G_{0,t+s}(u)=G_{t,t+s}(G_{0,t}(u)),\]
 and $G_{t,t+s}=\theta_{t}G_{0,s},$ then \[
a_{1}^{t+s}(u,\tau_{1})=T_{s}a_{1}^{t}(u,\tau_{1})1_{\tau_{1}<t}+\]
\[
a_{1}^{s}(T_{t}u,\tau_{1}-t)1_{t\leq\tau_{1}\leq t+s}.\qquad(1.4)\]
 Using this relation one can get \[
\begin{split} & a_{1}^{t}(u,0)=T_{t-\tau_{1}}a_{1}^{\tau_{1}}(u,0),\\
 & a_{1}^{t}(u,\tau_{1})=a_{1}^{t-\tau_{1}}(T_{\tau_{1}}u,0).\end{split}
\qquad(1.5)\]
 The condition of the theorem imply that for $v=T_{\tau_{1}}u$ and
every $h\in H_{0}$ there exists the limit \[
A(v)h=\lim_{t\to0+}\frac{1}{t}EG_{0,t}(v)(\xi,h\otimes1_{[0;\: t]})=\]
 \[
=\lim_{t\to0+}\frac{1}{t}\int_{0}^{t}a_{1}^{t}(v,\tau_{1})hd\tau_{1}.\]
 Now, by continuity of $a_{1},$ \[
a_{1}^{0}(T_{\tau_{1}}u,0)=A(T_{\tau_{1}}u).\]
 Finally, \[
a_{1}^{t}(u,\tau_{1})=T_{t-\tau_{1}}AT_{\tau_{1}}u.\]
 The case $k\geq2$ can be proved by induction. Suppose, that we have
the representation (1.3) for $a_{j}^{t},\: j\leq k.$ Consider $a_{k+1}^{t+s}.$
Using the multiplicative and homogeneity properties one can get \[
a_{k+1}^{t+s}(u,\tau_{1},\ldots,\tau_{k+1})1_{\{0\leq\tau_{1}\leq\ldots\leq\tau_{k}\leq t\leq\tau_{k+1}\leq t+s\}}=\]
 \[
=a_{1}^{s}(a_{k}^{t}(u,\tau_{1},\ldots,\tau_{k}),\tau_{k+1}-t)=\]
 \[
=T_{s+t-\tau_{k+1}}AT_{\tau_{k+1}-t}a_{k}^{t}(u,\tau_{1},\ldots,\tau_{k})=\]
 \[
=T_{s+t-\tau_{k+1}}AT_{\tau_{k+1}-t}T_{t-\tau_{k}}A\ldots AT_{\tau_{1}}u=\]
 \[
=T_{s+t-\tau_{k+1}}AT_{\tau_{k+1}-\tau_{k}}A\ldots AT_{\tau_{1}}u.\]
 The theorem is proved. }{\large \par}

{\large Consider some examples of application of the representation
(1.3).}{\large \par}

\textbf{\large Example 1.6.}{\large{} Consider the multiplicative functional
from Example 1.3. Suppose that the coefficients $a,\; b$ have infinitely
many derivatives. Now it can be proved, that $x(u,t)$ has infinitely
many stochastic derivatives {[}8{]}. Consequently for a smooth function
$f$ the first kernel in the Itô-Wiener expansion of $f(x(u,t))$
can be expressed as follows \[
a_{1}^{t}(u,\tau)=EDf(x(u,t))(\tau).\qquad(1.6)\]
 Indeed, for an arbitrary $h\in L_{2}([0;+\infty))$ \[
\int_{0}^{t}a_{1}^{t}(u,\tau)h(\tau)d\tau=Ef(x(u,t))\int_{0}^{t}h(\tau)dw(\tau)=\]
 \[
=E\int_{0}^{t}Df(x(u,t))(\tau)h(\tau)d\tau,\]
 which gives us the expression (1.6). The required continuity of $a_{1}$
follows from a well-known expression for the stochastic derivative
of $x$ {[}8{]}. As it was mentioned in Example 1.5, the operator
$A$ coincides with $b\frac{d}{du}$ on smooth functions. Finally,
the expression (1.3) turns into the well-known Krylov--Veretennikov
expansion {[}1{]} for $f(x(u,t))$ \[
f(x(u,t))=T_{t}f(u)+\sum_{k=1}^{\infty}\int_{\Delta_{k}(0;t)}T_{\tau_{1}}b\partial T_{\tau_{2}-\tau_{1}}\ldots\]
 \[
\ldots b\partial T_{t-\tau_{k}}f(u)dw(\tau_{1})\ldots dw(\tau_{k}).\]
 }\textbf{\textit{\large Remark}}\textbf{\large .}{\large{} The expression
(1.3) can be applied to multiplicative functionals, which are not
generated by a stochastic flow.}{\large \par}

\textbf{\large Example 1.7.}{\large{} Let $\mathcal{L}$ be a matrix
Lie group with the corresponding Lie algebra $\mathcal{A}$ with $\dim\mathcal{A}=n.$
Consider an $\mathcal{L}$-valued homogeneous multiplicative functional
$\{G_{s,t}\}$ from $\xi.$ Suppose that $\{G_{0,t}\}$ is a semimartingale
with respect to the filtration generated by $\xi.$ Let $\{G_{s,t}\}$
be continuous with respect to $s,t$ with probability one. It means,
in particular, that $\{G_{0,t}\}$ is a multiplicative Brownian motion
in $\mathcal{L}$ {[}9{]}. Then $G_{0,t}$ is a solution to the following
SDE \[
\begin{aligned} & dG_{0,t}=G_{0,t}dM_{t},\\
 & G_{0,0}=I.\end{aligned}
\]
 Here $\{M_{t};t\geq1\}$ is an $\mathcal{A}$-valued Brownian motion
obtained from $G$ by the rule {[}9{]} \[
M_{t}=P\mbox{-}\lim_{\Delta\to0+}\sum_{k=0}^{\left[\frac{t}{\Delta}\right]}(G_{k\Delta,(k+1)\Delta}-I).\qquad(1.7)\]
 Since $G_{0,t}$ is a semimartingale with respect to the filtration
of $\xi,$ then $M_{t}$ also has the same property. The representation
(1.7) shows that $M_{t}-M_{s}$ is measurable with respect to the
$\sigma$-field $\mathcal{F}_{s,t}$ and for arbitrary $r\geq0$ \[
\theta_{r}(M_{t}-M_{s})=M_{t+r}-M_{s+r}.\]
 Considering the Itô--Wiener expansion of $M_{t}-M_{s}$ one can easily
check, that \[
M_{t}=\int_{0}^{t}Zd\xi(\tau)\qquad(1.8)\]
 with a deterministic matrix $Z.$ We will prove (1.8) for the one-dimensional
case. Suppose that $M_{t}$ has the following Itô-Wiener expansion
with respect to $\xi$}{\large \par}

\[
M_{t}=\sum_{k=1}^{\infty}\int_{\Delta_{k}(t)}a_{k}(t,\tau_{1},\ldots,\tau_{k})d\xi(\tau_{1})\ldots d\xi(\tau_{k}).\]
 {\large Then for $k\geq2$ the corresponding kernel $a_{k}$ satisfies
relation\[
a_{k}(t+s,\tau_{1},\ldots,\tau_{k})=a_{k}(t,\tau_{1},\ldots,\tau_{k})1_{\{\tau_{1},\ldots,\tau_{k}\leq t\}}+\]
\[
a_{k}(s,\tau_{1}-t,\ldots,\tau_{k}-t)1_{\{\tau_{1},\ldots,\tau_{k}\geq t\}}.\]
 Iterating this relation for $t=\sum_{j=1}^{n}\frac{t}{n}$ one can
verify that $a_{k}\equiv0.$ For $k=1$ the same arguments give $a_{k}\equiv const.$}{\large \par}

{\large Consequently, the equation for $G$ can be rewritten using
$\xi$ as \[
dG_{0,t}=G_{0,t}Zd\xi(t).\qquad(1.9)\]
 Now the elements of the Itô--Wiener expansion from Theorem 1.1 can
be determined as follows \[
T_{t}=EG_{0,t},\ A=Z.\]
 Consequently, \[
G_{0,t}=T_{t}+\sum_{k=1}^{\infty}\int_{\Delta_{k}(0;t)}T_{t-\tau_{k}}ZT_{\tau_{k}-\tau_{k-1}}\ldots ZT_{\tau_{1}}d\xi(\tau_{1})\ldots d\xi(\tau_{k}).\]
 }{\large \par}

\section{The Arratia flow}

{\large When trying to obtain an analog of the representation (1.3)
for a stochastic flow which is not generated by a stochastic differential
equation with smooth coefficients, we are faced with the difficulty
that there is no such a Gaussian random vector field, which would
generate the flow. This circumstance arise from the possibility of
coalescence of particles in the flow. We will consider one of the
best known examples of such stochastic flows, the Arratia flow. Let
us start with the precise definition.}{\large \par}

\textbf{\large Definition 2.1.}{\large{} The Arratia flow is a random
field $\{x(u,t);\: u\in\mathbb{R},\: t\geq0\},$ which has the properties}{\large \par}

{\large 1) all $x(u,\cdot),\: u\in\mathbb{R}$ are the Wiener martingales
with respect to the join filtration,}{\large \par}

{\large 2) $x(u,0)=u,\ u\in\mathbb{R},$}{\large \par}

{\large 3) for all $u_{1}\leq u_{2},\ t\geq0$ \[
x(u_{1},t)\leq x(u_{2},t),\]
 }{\large \par}

{\large 4) the joint characteristics equals \[
<x(u_{1},\cdot),\: x(u_{2},\cdot)>_{t}=\int_{0}^{t}1_{\{\tau(u_{1},u_{2})\leq s\}}ds,\]
 where \[
\tau(u_{1},u_{2})=\inf\{t:\: x(u_{1},t)=x(u_{2},t)\}.\]
 }{\large \par}

{\large It follows from the properties 1)--3), that individual particles
in the Arratia flow move as Brownian particles and coalesce after
meeting. Property 4) reflects the independence of the particles before
meeting. It was proved in {[}10{]}, that the Arratia flow has a modification,
which is a càdlàg process on $\mathbb{R}$ with the values in $C([0;+\infty)]).$
From now on, we assume that we are dealing with such a modification.
We will construct the Arratia flow using a sequence of independent
Wiener processes $\{w_{k}:\: k\geq1\}.$ Suppose that $\{r_{k};\: k\geq1\}$
are rational numbers on $\mathbb{R}.$ To construct the Arratia flow
put $w_{k}(0)=r_{k},\: k\geq1$ and define \[
x(r_{1},t)=w_{1}(t),\ t\geq0.\]
 If $x(r_{1},\cdot),\ldots,\: x(r_{n},\cdot)$ have already been constructed,
then define \[
\sigma_{n+1}=\inf\{t:\:\prod_{k=1}^{n}(x(r_{k},t)-w_{n+1}(t))=0\},\]
 \[
x(r_{n+1},t)=\begin{cases}
w_{n+1}(t),\ t\leq\sigma_{n+1}\\
x(r_{k^{*}},t),\ t\geq\sigma_{n+1},\end{cases}\]
 where \[
w_{n+1}(\sigma_{n+1})=x(r_{k^{*}},\sigma_{n+1}),\]
\[
k=\min\{l:\; w_{n+1}(\sigma_{n+1})=x(r_{l},\sigma_{n+1})\}.\]
 In this way we construct a family of the processes $x(r,\cdot),\: r\in\mathbb{Q}$
which satisfies conditions 1)--4) from Definition 2.1. }{\large \par}

\textbf{\large Lemma 2.1.}{\large{} For every $u\in\mathbb{R}$ the
random functions $x(r,\cdot)$ uniformly converge on compacts with
probability one as $r\to u.$ For rational $u$ the limit coincides
with $x(u,\cdot)$ defined above. The resulting random field $\{x(u,t);\: u\in\mathbb{R},\: t\geq0\}$
satisfies the conditions of Definition 2.1 }{\large \par}

\textbf{\large Proof.}{\large{} Consider a sequence of rational numbers
$\{r_{n_{k}};\: k\geq1\}$ which converges to some $u\in\mathbb{R}\setminus\mathbb{Q}.$
Without loss of generality one can suppose that this sequence decreases.
For every $t\geq0$, $\{x(r_{n_{k}},t);\: k\geq1\}$ converges with
probability one as a bounded monotone sequence. Denote \[
x(u,t)=\lim_{k\to\infty}x(r_{n_{k}},t).\]
 Note, that for arbitrary $r',\: r''\in\mathbb{Q}$ and $t\geq0$\[
E\sup_{[0;\ t]}(x(r',s)-x(r'',s))^{2}\leq C\cdot(|r'-r''|+(r'-r'')^{2}).\qquad(2.1)\]
 Here the constant $C$ does not depend on $r'$ and $r''.$ Inequality
(2.1) follows from the fact, that the difference $x(r',\cdot)-x(r'',\cdot)$
is a Wiener process with variance 2, started at $r'-r''$ and stopped
at 0. Monotonicity and (2.1) imply that the first assertion of the
lemma holds. Note, that for every $t\geq0$ \[
\mathcal{F}_{t}=\sigma(x(r,s);\; r\in\mathbb{Q},\: s\in[0;t])=\]
 \[
=\sigma(x(r,s);\: r\in\mathbb{R},\: s\in[0;t]).\]
 Using standard arguments one can easily verify, that for every $u\in\mathbb{R}$,
$x(u,\cdot)$ is a Wiener martingale with respect to the flow $(\mathcal{F}_{t})_{t\geq0},$
and that the inequality \[
x(u_{1},t)\leq x(u_{2},t)\]
 remains to be true for all $u_{1}\leq u_{2}.$ Consequently, for
all $u_{1},u_{2}\in\mathbb{R}$, $x(u_{1},\cdot)$ and $x(u_{2},\cdot)$
coincide after meeting. It follows from (2.1) and property 4) for
$x(r,\cdot)$ with rational $r,$ that \[
<x(u_{1},\cdot),x(u_{2},\cdot)>_{t}=0\]
 for \[
t<\inf\{s:\: x(u_{1},s)=x(u_{2},s)\}.\]
 Hence, the family $\{x(u,t);\: u\in\mathbb{R},\: t\geq0\}$ satisfies
Definition 2.1.}{\large \par}

{\large This lemma shows that the Arratia flow is generated by the
initial countable system of independent Wiener processes $\{w_{k};k\geq1\}.$
From this lemma one can easily obtain the following statement. }{\large \par}

\textbf{\large Corollary 2.1.}{\large{} The $\sigma$-field \[
\mathcal{F}_{0+}^{x}:=\bigcap_{t>0}\sigma(x(u,s);\: u\in\mathbb{R},\:0\leq s\leq t)\]
 is trivial modulo $P.$ }{\large \par}

{\large The proof of this statement follows directly from the fact
that the Wiener process has the same property {[}11{]}. }{\large \par}

\section{The Krylov--Veretennikov expansion for the $n$-point motion of the
Arratia flow}

{\large We begin this section with an analog of the Krylov--Veretennikov
expansion for the Wiener process stopped at zero. For the Wiener process
$w$ define the moment of the first hitting zero\[
\tau=\inf\{t:\; w(t)=0\}\]
 and put $\widetilde{w}(t)=w(\tau\wedge t)$. For a measurable bounded
$f:\mathbb{R}\to\mathbb{R}$ define \[
\widetilde{T}_{t}(f)(u)=E_{u}f(\widetilde{w}(t)).\]
 The following statement holds.}{\large \par}

\textbf{\large Lemma 3.1. }{\large For a measurable bounded function
$f:\mathbb{R}\to\mathbb{R}$ and $u\geq0$\[
f(\widetilde{w}(t))=\widetilde{T}_{t}f(u)+\]
 \[
+\sum_{k=1}^{\infty}\int_{\Delta_{k}(t)}\widetilde{T}_{t-r_{k}}\frac{\partial}{\partial v_{k}}\widetilde{T}_{r_{k}-r_{k-1}}\ldots\]
\[
\frac{\partial}{\partial v_{1}}\widetilde{T}_{r_{1}}f(v_{1})dw(r_{1})\ldots dw(r_{k}).\qquad(3.1)\]
 }{\large \par}

\textbf{\large Proof. }{\large Let} {\large us use the Fourier --
Wiener transform. Define for $\varphi\in C([0;+\infty),\;\mathbb{R})\bigcap L_{2}([0;+\infty),\;\mathbb{R})$
the stochastic exponent \[
\mathit{\mathcal{E}}(\varphi)=\exp\{\int_{0}^{+\infty}\varphi(s)dw(s)-\frac{1}{2}\int_{0}^{+\infty}\varphi(s)^{2}ds\}.\]
 Suppose that a random variable $\alpha$ has the Itô--Wiener expansion
\[
\alpha=a_{0}+\sum_{k=1}^{\infty}\int_{\Delta_{k}(t)}a_{j}(r_{1},\ldots,r_{k})dw(r_{1})\ldots dw(r_{k}).\]
 Then \[
E\alpha\mathit{\mathcal{E}}(\varphi)=a_{0}+\]
 \[
+\sum_{k=1}^{\infty}\int_{\Delta_{k}(t)}a_{k}(r_{1},\ldots,r_{k})\varphi(r_{1})\ldots\varphi(r_{k})dr_{1}\ldots dr_{k}.\qquad(3.2)\]
 Consequently, to find the Itô-Wiener expansion of $\alpha$ it is
enough to find $E\alpha\mathit{\mathcal{E}}(\varphi)$ as an analytic
functional from $\varphi.$ Note that \[
E_{u}f(\widetilde{w}(t))\mathit{\mathcal{E}}(\varphi)=E_{u}f(\widetilde{y}(t)),\]
 where the process $\widetilde{y}$ is obtained from the process \[
y(t)=w(t)+\int_{0}^{t}\varphi(r)dr\]
 in the same way as $\widetilde{w}$ from $w.$ To find $E_{u}f(\widetilde{y}(t))$
consider the case when $f$ is continuous bounded function with $f(0)=0.$
Let $F$ be the solution to the following boundary problem on $[0;\;+\infty)\times[0;\; T]$
\[
\frac{\partial}{\partial t}F(u,t)=-\frac{1}{2}\frac{\partial^{2}}{\partial u^{2}}F(u,t)-\varphi(t)\frac{\partial}{\partial u}F(u,t),\qquad(3.3)\]
 \[
F(u,T)=f(u),\; F(0,s)=0,\; s\in[0;\; T],\]
 \[
F\in C^{2}((0;\;+\infty)\times(0;\; T))\cap C([0;\;+\infty)\times[0;\; T]).\]
 Then $F(u,0)=E_{u}f(\widetilde{y}(T)).$ To check this relation note,
that $F$ satisfies the relation\[
\frac{\partial}{\partial u}F(0,s)=\frac{\partial^{2}}{\partial u^{2}}F(0,s)=0,\; s\in[0;\; T].\]
 Consider the process $F(\widetilde{y}(s),s)$ on the interval $[0;\; T].$
Using Itô's formula one can get\[
F(\widetilde{y}(T),T)=F(u,0)+\int_{0}^{T\wedge\tau}(\frac{1}{2}\frac{\partial^{2}}{\partial u^{2}}F(\widetilde{y}(s),s)+\]
 \[
\varphi(s)\frac{\partial}{\partial u}F(\widetilde{y}(s),s))--(\frac{1}{2}\frac{\partial^{2}}{\partial u^{2}}F(\widetilde{y}(s),s)+\varphi(s)\frac{\partial}{\partial u}F(\widetilde{y}(s),s))ds+\]
 \[
+\int_{0}^{T\wedge\tau}\frac{\partial}{\partial u}F(\widetilde{y}(s),s)dw(s).\]
 Consequently\[
F(u,0)=E_{u}f(\widetilde{y}(T)).\]
 The problem (3.3) can be solved using the semigroup $\{\widetilde{T}_{t};\; t\geq0\}.$
It can be obtained from (3.3) that\[
F(u,s)=\widetilde{T}_{T-s}f(u)+\int_{s}^{T}\varphi(r)\widetilde{T}_{r-s}\frac{\partial}{\partial u}F(u,r)dr.\qquad(3.4)\]
 Solving (3.4) by the iteration method one can get the series\[
F(u,s)=\widetilde{T}_{T-s}f(u)+\]
 \[
+\sum_{k=1}^{\infty}\int_{\Delta_{k}(s;\; T)}\widetilde{T}_{r_{1}-s}\frac{\partial}{\partial v_{1}}\widetilde{T}_{r_{2}-r_{1}}\ldots\]
\[
\frac{\partial}{\partial v_{k}}\widetilde{T}_{T-r_{k}}f(v_{k})\varphi(r_{1})\ldots\varphi(r_{k})dr_{1}\ldots dr_{k}.\]
 The last formula means that the Itô-Wiener expansion of $f(\widetilde{w}(t))$
has the form\[
f(\widetilde{w}(t))=\widetilde{T}_{t}f(u)+\sum_{k=1}^{\infty}\int_{\Delta_{k}(t)}\widetilde{T}_{r_{1}}\frac{\partial}{\partial v_{1}}\widetilde{T}_{r_{2}-r_{1}}\ldots\]
\[
\frac{\partial}{\partial v_{k}}\widetilde{T}_{t-r_{k}}f(v_{k})dw(r_{1})\ldots dw(r_{k}).\qquad(3.5)\]
 To consider the general case note that for $t>0$ and $c\in\mathbb{R}$\[
\frac{\partial}{\partial v}\widetilde{T}_{t}c\equiv0.\]
 Consequently (3.5) remains to be true for an arbitrary bounded continuous
$f.$ Now the statement of the lemma can be obtained using the approximation
arguments. The lemma is proved. }{\large \par}

{\large The same idea can be used to obtain the Itô-Wiener expansion
for a function from the Arratia flow. The $n-$point motion of the
Arratia flow was constructed in Section 2 from independent Wiener
processes. Consequently, a function from this $n-$point motion must
have the Itô-Wiener expansion in terms of these processes. We will
treat such expansion as the Krylov-Veretennikov expansion for the
Arratia flow.}{\large \par}

{\large Here there is a new circumstance compared to the case when
the flow is generated by SDE with smooth coefficients. Namely, there
are many different ways to construct the trajectories of the Arratia
flow from the initial Wiener processes, and the form of the Itô-Wiener
expansion will depend on the way of constructing the trajectories.
In {[}12{]} Arratia described different ways of constructing the colliding
Brownian motions from independent Wiener processes. We present here
a more general approach by considering a broad class of constructions,
and find the Itô-Wiener expansion for it. To describe our method we
will need some preliminary notations and definitions. }{\large \par}

\textbf{\large Definition 3.1. }{\large An arbitrary set of the kind
$\{i,i+1,\ldots,j\}$, where $i,j\in\mathbb{N}$,$i\leq j$ is called
a block.}{\large \par}

\textbf{\large Definition 3.2.}{\large{} A representation of the block
$\{1,2,\ldots,n\}$ as a union of disjoint blocks is called a partition
of the block $\{1,2,\ldots,n\}$. }{\large \par}

\textbf{\large Definition 3.3. }{\large We say that a partition $\pi_{2}$
follows from a partition $\pi_{1}$ if it coincides with $\pi_{1}$
or if it is obtained by the union of two subsequent blocks from $\pi_{1}$. }{\large \par}

{\large We will consider a sequences of partitions $\{\pi_{0},\ldots,\pi_{l}\}$
where $\pi_{0}$ is a trivial partition, $\pi_{0}=\{\{1\},\{2\},\ldots,\{n\}\}$
and every $\pi_{i+1}$ follows from $\pi_{i}.$ The set of all such
sequences will be denoted by $R$. Denote by $R_{k}$ the set of all
sequences from $R$ that have exactly $k$ matching pairs: $\pi_{i}=\pi_{i+1}$.
The set $R_{0}$ of strongly decreasing sequences we denote by $\breve{R}.$
For every sequence $\{\pi_{0},\ldots,\pi_{k}\}$ from $\breve{R}$
each $\pi_{i+1}$ is obtained from $\pi_{i}$ by the union of two
subsequent blocks. It is evident, that the length of every sequence
from $\breve{R}$ is less or equal to $n.$ Let us associate with
every partition $\pi$ a vector $\vec{\lambda_{\pi}}\in\mathbb{R}^{n}$
with the next property. For each block $\{s,\ldots,t\}$ from $\pi$
the following relation holds \[
\sum_{q=s}^{t}\lambda_{\pi q}^{2}=1.\]
 We will use the mapping $\vec{\lambda}$ as a rule of constructing
the $n-$point motion of the Arratia flow. Suppose now, that $\{w_{k};\: k=1,\ldots,n\}$
are independent Wiener processes starting at the points $u_{1}<\ldots<u_{n}.$
We are going to construct the trajectories $\{x_{1},\ldots,x_{n}\}$
of the Arratia flow starting at $u_{1}<\ldots<u_{n}$ from the pieces
of the trajectories of $\{w_{k};\: k=1,\ldots,n\}$. Assume that we
have already built the trajectories of $\{x_{1},\ldots,x_{n}\}$ up
to a certain moment of coalescence $\tau$. At this moment a partition
$\pi$ of $\{1,2,\ldots,n\}$ naturally arise. Two numbers $i$ and
$j$ belong to the same block in $\pi$ if and only if $x_{i}(\tau)=x_{j}(\tau).$
Consider one block $\{s,\ldots,t\}$ in $\pi$. Define the processes
$x_{s},\ldots,x_{t}$ after the moment $\tau$ and up to the next
moment of coalescence in the whole system $\{x_{1},\ldots,x_{n}\}$
by the rule\[
x_{i}(t)=x_{i}(\tau)+\sum_{q=s}^{t}\lambda_{\pi q}(w_{q}(t)-w_{q}(\tau)).\]
 Proceeding in the same way, we obtain the family $\{x_{k},\; k=1,\ldots,n\}$
of continuous square integrable martingales with respect to the initial
filtration, generated by $\{w_{k};\: k=1,\ldots,n\}$ with the following
properties:}{\large \par}

{\large 1) for} {\large every $k=1,\ldots,n,\quad x_{k}(0)=u_{k},$}{\large \par}

{\large 2) for every $k=1,\ldots,n-1,\quad x_{k}(t)\leq x_{k+1}(t)$,}{\large \par}

{\large 3) the joint characteristic of $x_{i}$ and $x_{j}$ satisfies
relation \[
d<x_{i},\; x_{j}>(t)=\mathbf{1}_{t\geq\tau_{ij}},\]
 where $\tau_{ij}=\inf\{s:\; x_{i}(s)=x_{j}(s)\}.$}{\large \par}

{\large It can be proved {[}13{]} that the processes $\{x_{k},\; k=1,\ldots,n\}$
are the $n-$point motion of the Arratia flow starting from the points
$u_{1}<\ldots<u_{n}.$ We constructed it from the independent Wiener
processes $\{w_{k};\: k=1,\ldots,n\}$ and the way of construction
depends on the mapping $\vec{\lambda}.$ To describe the Itô-Wiener
expansion for functions from $\{x_{k}(t),\; k=1,\ldots,n\}$ it is
necessary to introduce operators related to a sequence of partitions
$\tilde{\pi}\in\breve{R}.$ Denote by $\tau_{0}=0<\tau_{1}<\ldots<\tau_{n-1}$
the moments of coalescence for $\{x_{k}(t),\; k=1,\ldots,n\}$ and
by $\tilde{\nu}=\{\pi_{0},\nu_{1},\ldots,\nu_{n-1}\}$ related random
sequence of partitions. Namely, the numbers $i$ and $j$ belong to
the same block in the partition $\nu_{k}$ if and only if $x_{i}(t)=x_{j}(t)$
for $\tau_{k}\leq t.$ Define for a bounded measurable function $f:\mathbb{R}^{n}\to\mathbb{R}$
\[
T_{t}^{\tilde{\pi}}f(u_{1},\ldots,u_{n})=Ef(x_{1}(t),\ldots,x_{n}(t))\mathbf{1}_{\{\nu_{1}=\pi_{1},\ldots,\nu_{k}=\pi_{k},\;\tau_{k}\leq t<\tau_{k+1}\}}.\]
}{\large \par}

{\large Now let $\kappa$ be an arbitrary partition and let $u_{1}\leq u_{2}\leq\ldots\leq u_{n}$
be such, that $u_{i}=u_{j}$ if and only if $i$ and $j$ belong to
the same block in $\kappa.$ One can define formally the $n-$point
motion of the Arratia flow starting at $u_{1}\leq u_{2}\leq\ldots\leq u_{n}$,
assuming that the trajectories that start at coinciding points, also
coincide. Then for the strongly decreasing sequence of partitions
$\tilde{\pi}=\{\kappa,\pi_{1},\ldots,\pi_{k}\}$ the operator $T_{t}^{\tilde{\pi}}$
is defined by the same formula as above. }{\large \par}

{\large The next theorem is the Krylov-Veretennikov expansion for
the $n-$point motion of the Arratia flow.}{\large \par}

\textbf{\large Theorem 3.1. }{\large For a bounded measurable function
$f:\mathbb{R}^{n}\to\mathbb{R}$ the following representation takes
place\[
f(x_{1}(t),\ldots,x_{n}(t))=\sum_{\tilde{\pi}\in\breve{R}}T_{t}^{\tilde{\pi}}f(u_{1},\ldots,u_{n})+\]
 \[
\sum_{i=1}^{n}\sum_{\tilde{\pi}\in R_{1}}\lambda_{\pi_{1}i}\int_{0}^{t}T_{s_{1}}^{\tilde{\pi}_{1}}\partial_{i}T_{t-s_{1}}^{\tilde{\pi}_{2}}f(u_{1},\ldots,u_{n})dw_{i}(s_{1})+\]
 \[
\sum_{i_{1},i_{2}=1}^{n}\sum_{\tilde{\pi}\in R_{2}}\lambda_{\pi_{1}i_{1}}\lambda_{\pi_{2}i_{2}}\int_{\triangle_{2}(t)}T_{s_{1}}^{\tilde{\pi}_{1}}\partial_{i_{1}}T_{s_{2}-s_{1}}^{\tilde{\pi}_{2}}\partial_{i_{2}}T_{t-s_{2}}^{\tilde{\pi}_{3}}\]
 \[
f(u_{1},\ldots,u_{n})dw_{i_{1}}(s_{1})dw_{i_{2}}(s_{2})+\]
 \[
\sum_{i_{1},...,i_{k}=1}^{n}\sum_{\tilde{\pi}\in R_{k}}\prod_{j=1}^{k}\lambda_{\pi_{j}i_{j}}\int_{\triangle_{k}(t)}T_{s_{1}}^{\tilde{\pi}_{1}}\partial_{i_{1}}T_{s_{2}-s_{1}}^{\tilde{\pi}_{2}}...\partial_{i_{k}}T_{t-s_{k}}^{\tilde{\pi}_{k+1}}\]
 \[
f(u_{1},\ldots,u_{n})dw_{i_{1}}(s_{1})...dw_{i_{k}}(s_{k})+...\]
 In this formula we use the following notations. For a sequence $\tilde{\pi}\in R_{k}$
partitions $\pi_{1},...,\pi_{k}$ are the left elements of equalities
from $\tilde{\pi}=\{...\pi_{1}=...\pi_{2}=...\pi_{k}=...\}$ and $\tilde{\pi}_{1},...,\tilde{\pi}_{k+1}$
are strictly decreasing pieces of $\tilde{\pi}.$ The symbol $\partial_{i}$
denotes differentiation with respect to a variable corresponding to
the block of partition, which contains $i.$ For example, if $i\in\{s,...,t\}$
then\[
\partial_{i}f=\sum_{q=s}^{q=t}f_{q}^{\prime}.\]
 }{\large \par}

{\large The proof of the theorem can be obtained by induction, adopting
ideas of Lemma 3.1. One has to consider subsequent boundary value
problems and then use the probabilistic interpretation of the Green's
functions for these problems. The corresponding routine calculations
are omitted. }{\large \par}

\textbf{\large References}{\large \par}

{\large 1. N. V. Krylov, A. Yu. Veretennikov. Explicit formulae for
the solutions of the stochastic differential equations. Math. USSR
Sb. 29 (1976), No. 2, pp. 239-256.}{\large \par}

{\large 2. A.V.Skorokhod. Random linear operators. D.Reidel Publishing
Company, 1983. -- Dordrecht, Holland. -- 198 p. }{\large \par}

{\large 3. H.Kunita. Stochastic flows and stochastic differential
equations. -- Cambridge University Press, 1990. -- 346 p. }{\large \par}

{\large 4. T. E. Harris. Coalescing and noncoalescing stochastic flows
in $\mathbb{R}_{1}$. Stochastic Processes and their Applications
17(1984), pp. 187 - 210.}{\large \par}

{\large 5. D.G. Aronson. Bounds for the fundamental solution of a
parabolic equation.- Bull. Amer. Math. Soc.-1967.-P. 890-896.}{\large \par}

{\large 6. A.A.Dorogovtsev. Stochastic analysis and random maps in
Hilbert space. -- Utrecht: VSP, 1994. -- 110 p. }{\large \par}

{\large 7. S.Janson. Gaussian Hilbert spaces. Cambridge University
Press, 1997. -- X+340 pp. }{\large \par}

{\large 8. S.Watanabe. Lectures on stochastic differential equations
and Malliavin calculus. -- Tata Institute of Fundamental Research,
Bombay, 1984. -- III+111 pp. }{\large \par}

{\large 9. A.S.Holevo. An analog of the Itô decomposition for multiplicative
processes with values in a Lie group. --Sankhya: The Indian Journal
of Statistics. -- 1991, Vol.53, Ser. A, Pt.2. -- P. 158--161. }{\large \par}

{\large 10. A.A.Dorogovtsev. Some remarks on a Wiener flow with coalescence.
--Ukrainian mathematical journal. --2005, Volume 57, Number 10, p.
1550-1558.}{\large \par}

{\large 11. O. Kallenberg. Foundations of modern probability. --Springer-Verlag,
1997.--VI+535 pp.}{\large \par}

{\large 12. R. Arratia. Coalescing Brownian motion on the line. PhD
thesis. University of Wisconsin - Madison, 1979.}{\large \par}

{\large 13. V.V. Konarovskii. On Infinite System of Diffusing Particles
with Coalescing. -- Theory Probab. Appl. 55, pp. 134-144. }{\large \par}

{\large 14.  Y. Le Jan, O.Raimond. Flows, coalescence and noise. --
Ann.Probab. 32 (2004). -- P. 1247-1315.} 
\end{document}